\documentclass[11pt,letterpaper]{article}
\usepackage{amssymb, amsthm, amsmath}
\usepackage{graphicx}
\usepackage{hyperref}
\usepackage{subfig}

\newtheorem{theorem}{Theorem}
\newtheorem{observation}{Observation}

\begin{document}
\title{Subdivision by bisectors is dense in the space of all triangles}
\author{Steve Butler\thanks{Department of Mathematics, UCLA,
Los Angeles, CA 90095 ({\tt butler@math.ucla.edu}).}~\thanks{This work was done with support of an NSF Mathematical Sciences Postdoctoral Fellowship.} \and Ron Graham\thanks{Department of Computer Science and Engineering, University of California, San Diego,
La Jolla, CA 92093 ({\tt graham@ucsd.edu}).}}
\date{\empty}
\maketitle

\begin{abstract}
Starting with any nondegenerate triangle we can use a well defined interior point of the triangle to subdivide it into six smaller triangles.  We can repeat this process with each new triangle, and continue doing so over and over.  We show that starting with any arbitrary triangle, the resulting set of triangles formed by this process  contains triangles arbitrarily close (up to similarity) any given triangle when the point that we use to subdivide is the incenter.  We also show that the smallest angle in a ``typical'' triangle after repeated subdivision for many generations does not have the smallest angle going to $0$.
\end{abstract}

\section{Introduction}

\begin{figure}[htb]
\centering
\includegraphics[scale=0.8]{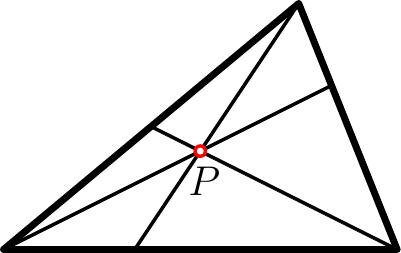}
\caption{How to subdivide the triangle given an interior point $P$.}
\label{fig:subdivide}
\end{figure}

Given a triangle and an interior point of the triangle we can divide the triangle up into six smaller triangles, also called daughters, by drawing line segments (or Cevians) from the vertex through the interior point to the opposite side (see Figure~\ref{fig:subdivide}).  The process can then be repeated with each new triangle with its own corresponding interior point and again repeated over and over.  When the interior point is the centroid this corresponds to barycentric subdivision.

Stakhovskii asked whether repeated barycentric subdivision for a starting triangle is dense is the space of all triangles, i.e., the space of triangles up to similarity where two triangles are $\epsilon$-close if the maximum difference between their corresponding angles is less than $\epsilon$.  This question was answered in the affirmative \cite{barany}.

\begin{theorem}[B\'ar\'any-Beardon-Carne]\label{thm:centroid}
Successive subdivisions of a non-degenerate triangle using the centroid point contains triangles which approximate arbitrarily closely (up to similarity) any given triangle.
\end{theorem}

Moreover, it was shown that almost all triangles became ``flat'' (in the sense that the largest angle approaches $\pi$ for most triangles as the number of subdivisions increase).  Ordin \cite{ordin} extended their results and showed that this also holds if we choose the interior point of the triangle to be $p_0 A+p_1 B+p_2 C$ where $A,B,C$ are the vertices with $p_0+p_1+p_2=1$ and $p_i>0$ (the centroid corresponds to $p_0=p_1=p_2={1\over3}$).

Our main result is to show that a similar statement holds if we choose the interior point to be the incenter, which can be found by taking the intersection of the angle bisectors.

\begin{theorem}\label{thm:incenter}
Succesive subdivisions of a non-degenerate triangle using the incenter point contains triangles which approximate arbitrarily closely (up to similarity) any given triangle.
\end{theorem}

We will see however that there is a difference in the behavior in that not almost all triangles become flat as in the centroid case.

We will proceed as follows.  In Section~\ref{sec:centroid} we will give a quick sketch of Theorem~\ref{thm:centroid}, while in Section~\ref{sec:incenter} we will give a proof of Theorem~\ref{thm:incenter} and establish several properties about this subdivision.  Finally, in Section~\ref{sec:conclusion} we will give some concluding remarks.

%
%

\section{Subdividing using the centroid}\label{sec:centroid}
\begin{figure}[ht]
\centering
\includegraphics[scale=0.8]{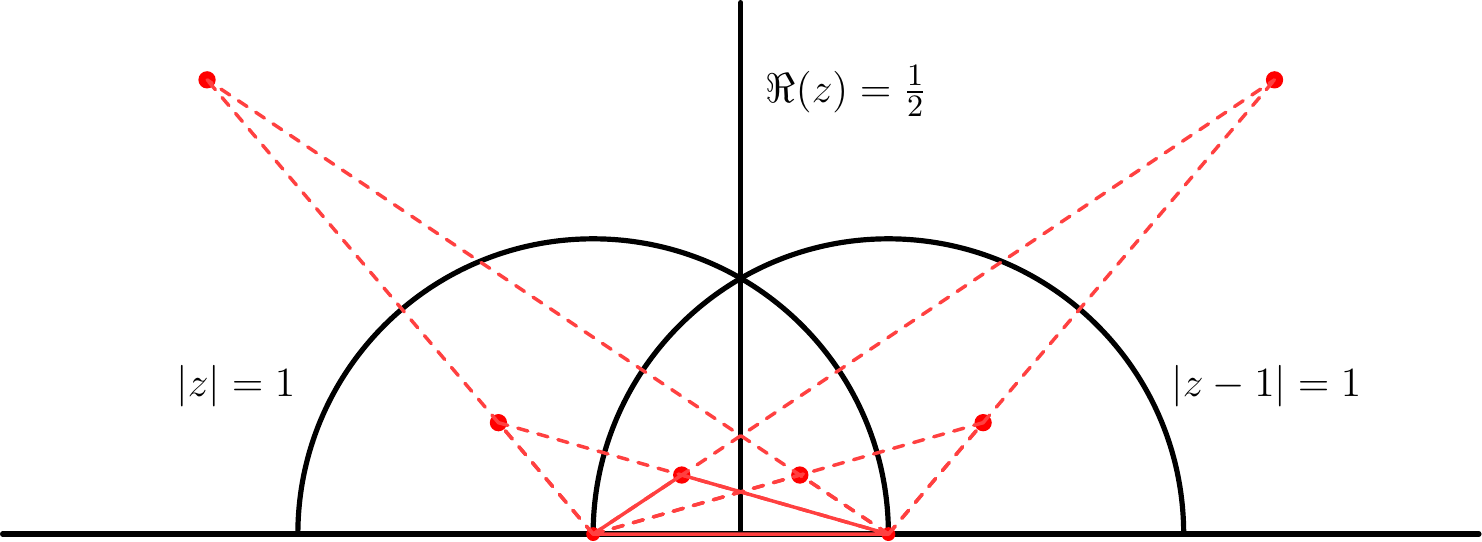}
\caption{Triangles as points in the hyperbolic plane.}
\label{fig:hyperbolic}
\end{figure}

In this section we give a quick sketch of the ideas behind Theorem~\ref{thm:centroid}.  The method of B\'ar\'any et al.\ \cite{barany} was to first associate triangles with points in  the hyperbolic half plane, namely each triangle $T$ is associated with (up to) six points $z$ in the hyperbolic upper half plane $\mathbb{H}$ as shown in Figure \ref{fig:hyperbolic}. This is done by placing {\it some}\/ edge of $T$ with vertices at $z=0$ and $z=1$ and the third vertex is located at the complex coordinate $z$ with positive imaginary part. Observe that reflecting $z$ across the three circles $\Re(z) = \frac{1}{2}$, $|z| = 1$, and $|z-1|=1$ induces a natural action of $S_3$ on the hyperbolic half plane $\mathbb{H}$ in which all six orientations of $T$ occur.

Now, the centroid point of a triangle with vertices at $0$, $1$ and $z$ is the point $(z+1)/3$, and so one of the corresponding daughters becomes $2(z+1)/3$ when normalized. The argument reduces to showing that the group of automorphisms of $\mathbb{H}$ generated by the map $B(z)=2(z+1)/3$ and the above $S_3$ action is dense in $Aut(\mathbb{H})$, in particular for any starting $z$ (i.e., any intitial triangle $T$) the set of all resulting points is dense in $\mathbb{H}$ (i.e., dense in the space of triangles).

Further, using results of Furstenberg \cite{furstenberg}, it follows that almost all random walks formed from products of $B(z)$ and elements of $S_3$ tend to infinity (in the hyperbolic plane) as the length of the product increases. This then implies that almost all of the $n$th generation daughters have smallest angle tending to $0$ as $n$ increases.  By different techniques, Robert Hough \cite{hough} was able to show that the largest angle approaches $\pi$ and moreover was able to give asymptotic bounds for the proportion of triangles with angles near $\pi$.

%
%

\section{Subdividing using the incenter}\label{sec:incenter}
\begin{figure}[ht]
\centering
\includegraphics[scale=0.8]{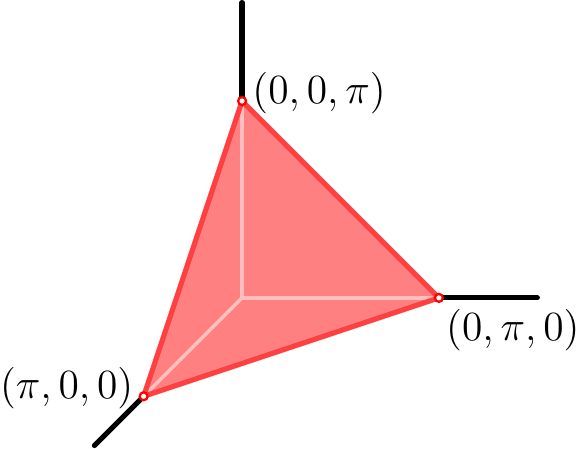}
\caption{Representing triangles as points on $\mathbf{P}$ in $\mathbb{R}^3$.}
\label{fig:H}
\end{figure}

The important step in the proof of Theorem~\ref{thm:centroid} was to find a way to associate triangles with points where the action of finding a daughter triangle was natural.  The first step in proving Theorem~\ref{thm:incenter} is to do the same.  However, we will find it more convenient to associate  each triangle with a point(s) in $\mathbb{R}^3$ where the coordinates are the angles.  The set of all possible triangles (including degenerate cases), denoted $\mathbf{P}$, is the intersection of the hyperplane $x+y+z=\pi$ with the first octant (see Figure~\ref{fig:H}).  Note that $\mathbf{P}$ also corresponds to a two dimensional equilateral triangle (see \cite{alexander,iterated,kingston,lax} for previous applications involving $\mathbf{P}$).

As noted in the introduction the incenter is found by the intersection of the angle bisectors.  So the angles of the new triangles created by subdivision are linear combinations of the angles of the original triangles (hence the reason it is more convenient to work with $\mathbf{P}$).  In particular, if we let $\mathbf{t}=(\alpha,\beta,\gamma)^*$ denote a triangle then the six new triangles are found by $M_i\mathbf{t}$ where
\begin{center}
$\displaystyle{
M_1=\begin{pmatrix}1/2&0&0\\1/2&1/2&0\\0&1/2&1\end{pmatrix},\qquad
M_2=\begin{pmatrix}1/2&1/2&0\\0&1/2&0\\1/2&0&1\end{pmatrix},\qquad
M_3=\begin{pmatrix}1&0&1/2\\0&1/2&0\\0&1/2&1/2\end{pmatrix},
}$\\
$\displaystyle{
M_4=\begin{pmatrix}1&1/2&0\\0&1/2&1/2\\0&0&1/2\end{pmatrix},\qquad
M_5=\begin{pmatrix}1/2&0&1/2\\1/2&1&0\\0&0&1/2\end{pmatrix},\qquad
M_6=\begin{pmatrix}1/2&0&0\\0&1&1/2\\1/2&0&1/2\end{pmatrix}.
}$
\end{center}

\begin{figure}[h]
\includegraphics[width=0.75in]{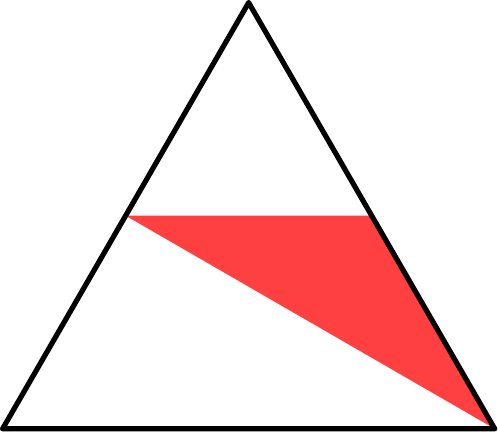}\hfil
\includegraphics[width=0.75in]{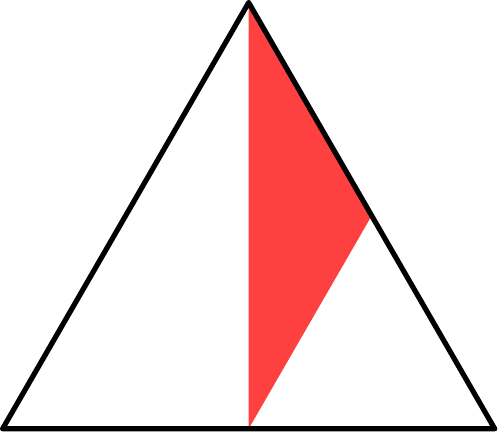}\hfil
\includegraphics[width=0.75in]{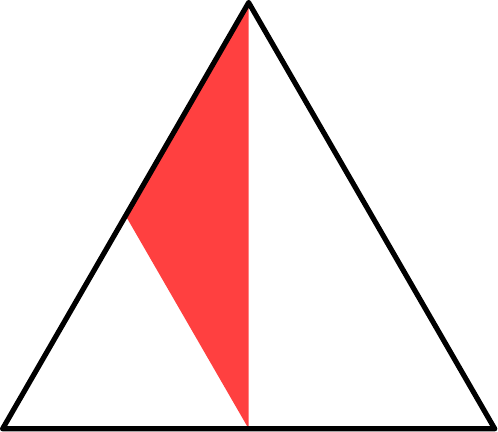}\hfil
\includegraphics[width=0.75in]{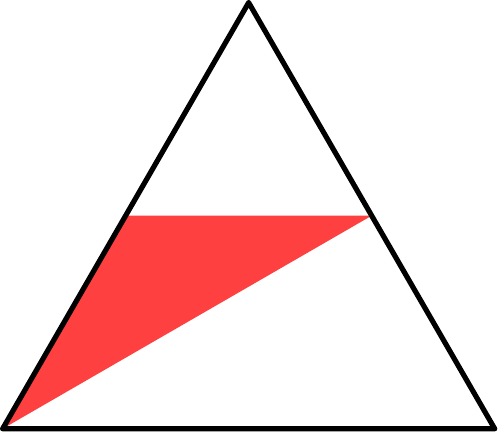}\hfil
\includegraphics[width=0.75in]{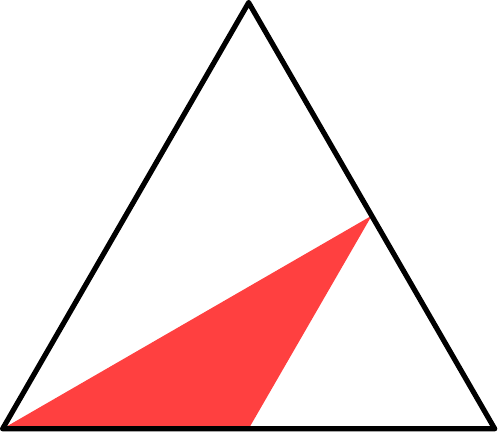}\hfil
\includegraphics[width=0.75in]{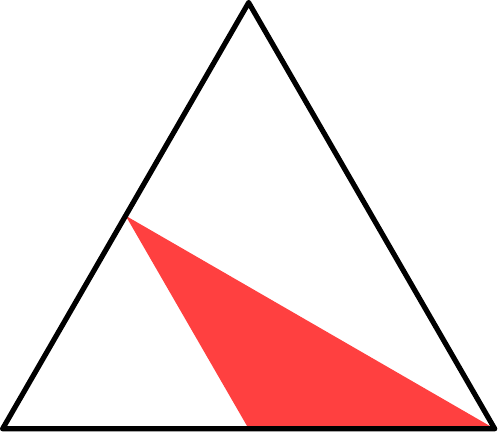}
\caption{The image of $\mathbf{P}$ under the six maps $M_i$.}
\label{fig:sixmaps}
\end{figure}

\begin{observation}\label{obs:cover}
The union of the $M_i\mathbf{P}$ covers $\mathbf{P}$.
\end{observation}

This can be seen by examining Figure~\ref{fig:sixmaps}.  Alternatively, this says that every point in $\mathbf{P}$ has a preimage in $\mathbf{P}$ under some $M_i$.  For the case that $\mathbf{t}=(\alpha,\beta,\gamma)^*$ with $\alpha\leq\beta\leq\gamma$ we have that $(M_1)^{-1}\mathbf{t}=(2\alpha,2\beta-2\alpha,\alpha-\beta+\gamma)^*$ is a preimage of $\mathbf{t}$ in $\mathbf{P}$. Other possible arrangements for the ordering of $\alpha,\beta,\gamma$ can be handled with the remaining $M_i$.

\begin{observation}\label{obs:shrink}
If we let $\|{\cdot}\|$ denote Euclidean distance, then $\|M_i(\mathbf{t}-\mathbf{s})\|\leq {\sqrt{3}\over 2}\|\mathbf{t}-\mathbf{s}\|.$
\end{observation}

To see this we can put $\mathbf{P}$ into $\mathbb{R}^2$ by $\mathbf{t}=(\alpha,\beta,\gamma)^*\to\big({\alpha+2\beta\over\sqrt3},\alpha\big)^*$; note that this will preserve distance.  Under this map we would also have that $M_1\mathbf{t}\to\big({3\alpha+2\beta\over2\sqrt3},{\alpha\over2}\big)^*$.  If we now let $\mathbf{s}=(\alpha',\beta',\gamma')$ then a calculation shows
\[
{3\over4}\|\mathbf{t}-\mathbf{s}\|^2-\|M_1\mathbf{t}-M_1\mathbf{s}\|^2={2\over3}\big(\beta-\beta'\big)^2\geq 0.
\]
The result now follows for $M_1$ and similar calculations establish it for the remaining $M_i$.

We now prove Theorem~\ref{thm:incenter}.  Let $\mathbf{q}$ be the initial triangle we apply subdivision to.  We need to show that for any triangle $\mathbf{t}$ and $\epsilon>0$ there is some sequence of $i_j$ so that
\[
\big\|(M_{i_1}M_{i_2}\cdots M_{i_k}\mathbf{q})-\mathbf{t}\big\|<\epsilon.
\]
Choose $k$ sufficiently large so that $\pi\big({\sqrt3\over2}\big)^k<\epsilon$.  By Observation~\ref{obs:cover} we can successively find a $k$th generation preimage of $\mathbf{t}$ in $\mathbf{P}$, which corresponds to multiplying by an appropriate $(M_i)^{-1}$ at each step.  Denote this preimage by
\[
(M_{i_k})^{-1}\cdots(M_{i_2})^{-1}(M_{i_1})^{-1}\mathbf{t},
\]
(where the $i_j$ are chosen according to how we construct the preimage).  Repeatedly using Observation~\ref{obs:shrink} we have
\begin{eqnarray*}
\big\|(M_{i_1}M_{i_2}\cdots M_{i_k}\mathbf{q})-\mathbf{t}\big\|&=&
\big\|M_{i_1}\big(M_{i_2}\cdots M_{i_k}\mathbf{q}-(M_{i_1})^{-1}\mathbf{t}\big)\big\|\\
&\leq&{\sqrt3\over2}\big\|M_{i_2}\cdots M_{i_k}\mathbf{q}-(M_{i_1})^{-1}\mathbf{t}\big\|\\
&\leq&\cdots\\
&\leq&\big({\sqrt{3}\over2}\big)^k\big\|\mathbf{q}-(M_{i_k})^{-1}\cdots(M_{i_2})^{-1}(M_{i_1})^{-1}\mathbf{t}\big\|\\
&\leq&\pi\big({\sqrt{3}\over2}\big)^k~<~\epsilon.
\end{eqnarray*}
In the last step we used that points in $\mathbf{P}$ are at most distance $\pi$ apart.  This finishes the proof of Theorem~\ref{thm:incenter}.

\subsection*{The limiting distribution}
In fact more can be said about the iterated subdivision of triangles using the incenter.  Namely, since the maps $M_i$ are contracting with Lipschitz constant $\sqrt3/2$ then it follows (see \cite{persi}) that there is a fixed limiting distribution on $\mathbf{P}$ that the process converges to.  Further it converges exponentially.  

\begin{figure}[ht]
\centering
\hfil\subfloat[]{\includegraphics[width=2.5in]{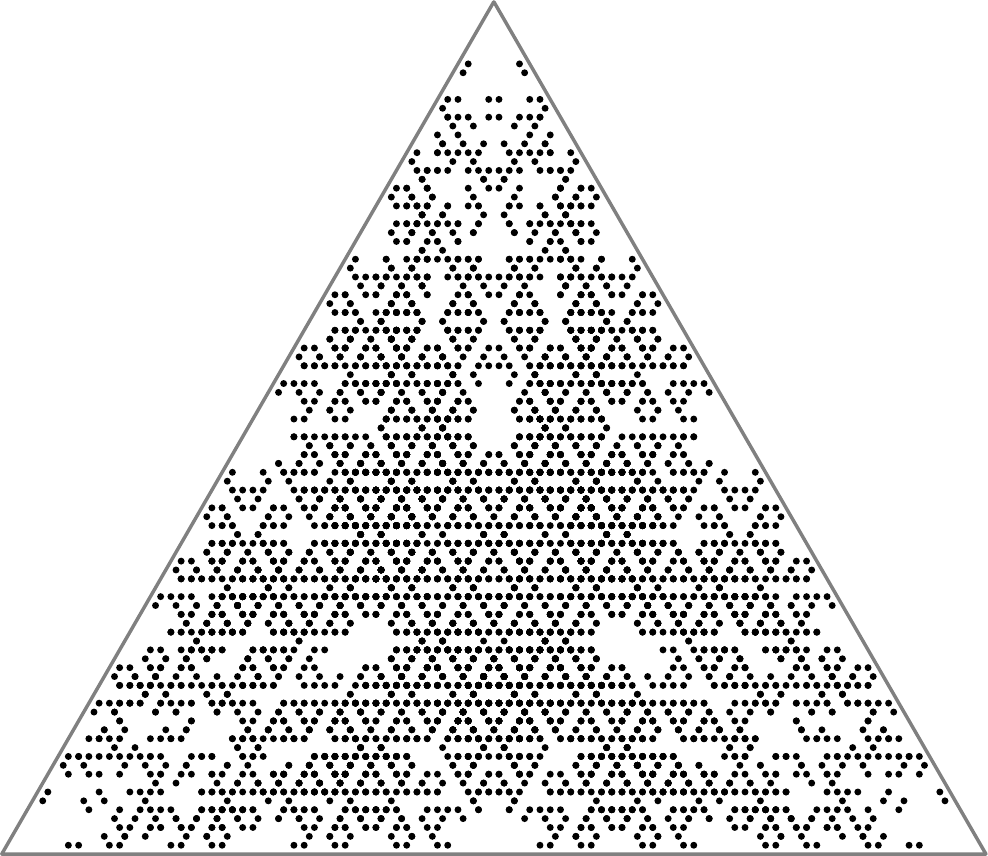}\label{fig:bb5}} \hfil
\hfil\subfloat[]{\includegraphics[width=2.5in]{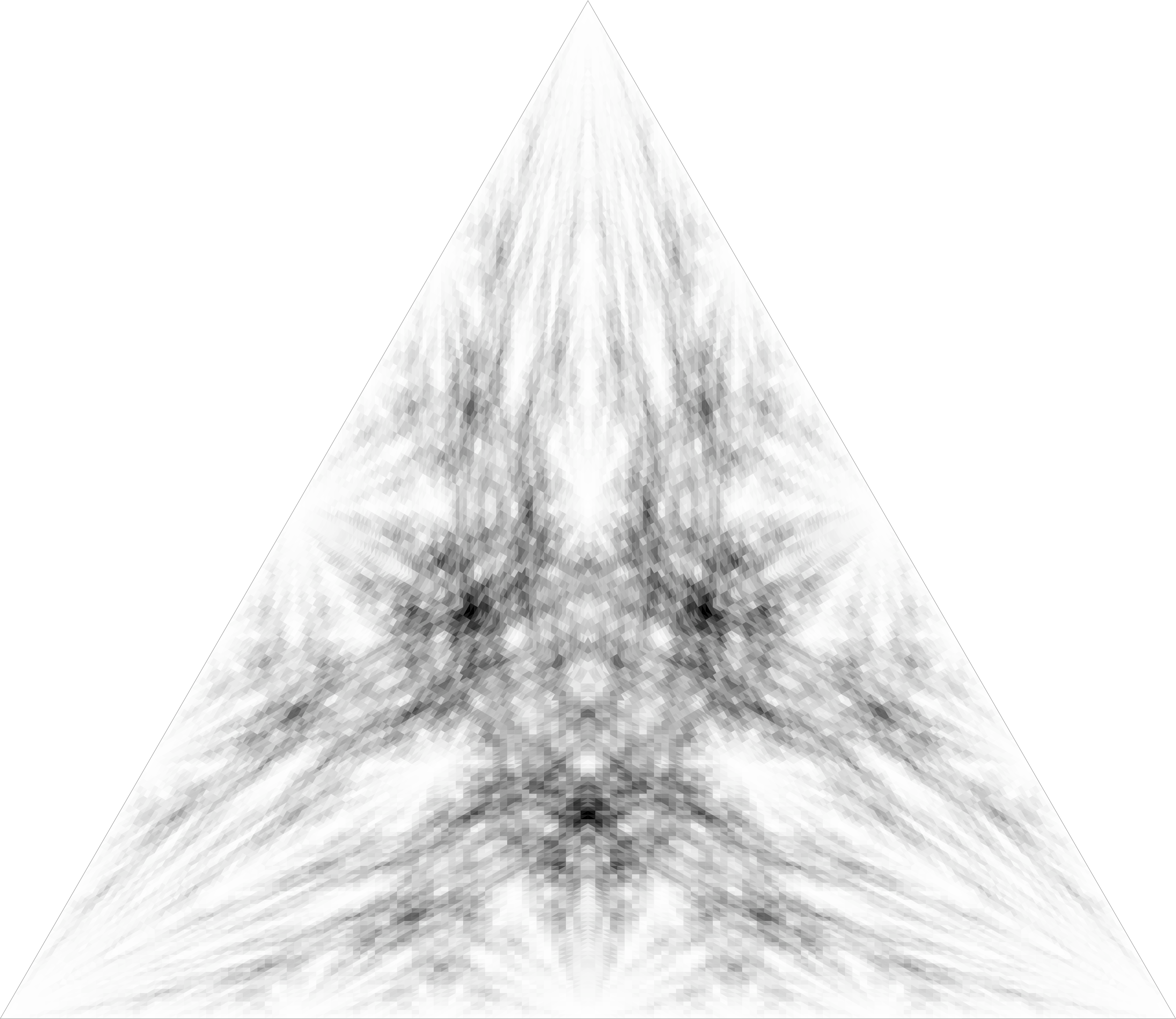}\label{fig:bb12}}\hfil
\caption{The distribution of the $n$th generation daughters using the incenter.}
\label{fig:bi_equ}
\end{figure}

To get some sense of what this limiting distribution looks like we can simply start with any triangle (in our case we will use an equilateral triangle) and plot all of the $n$th generation daughters for some $n$ in $\mathbf{P}$.  This is done for $n=5$ in Figure~\ref{fig:bb5}.

Examining Figure~\ref{fig:bb5} we see that the daughters seem to fill in most of $P$ (agreeing with Theorem~\ref{thm:incenter}). However, a patient count will reveal that there are far fewer than $6^5$ triangles in Figure~\ref{fig:bb5}. This is because some points have been mapped onto several times (a consequence of starting with such a symmetric triangle).  So to get a better sense of the limiting distribution instead of plotting the individual triangles in $\mathbf{P}$ it is better to look at a histogram.  We will divide $\mathbf{P}$ into a large number of small regions and then shade each region according to the number of triangles that fall into that region, the darker a region is the more triangles fall into that region.  In Figure~\ref{fig:bb12} we give the histogram for $n=12$ generations starting with the equilateral triangle.

\begin{figure}[ht]
\centering
\includegraphics[scale=0.4,angle=270]{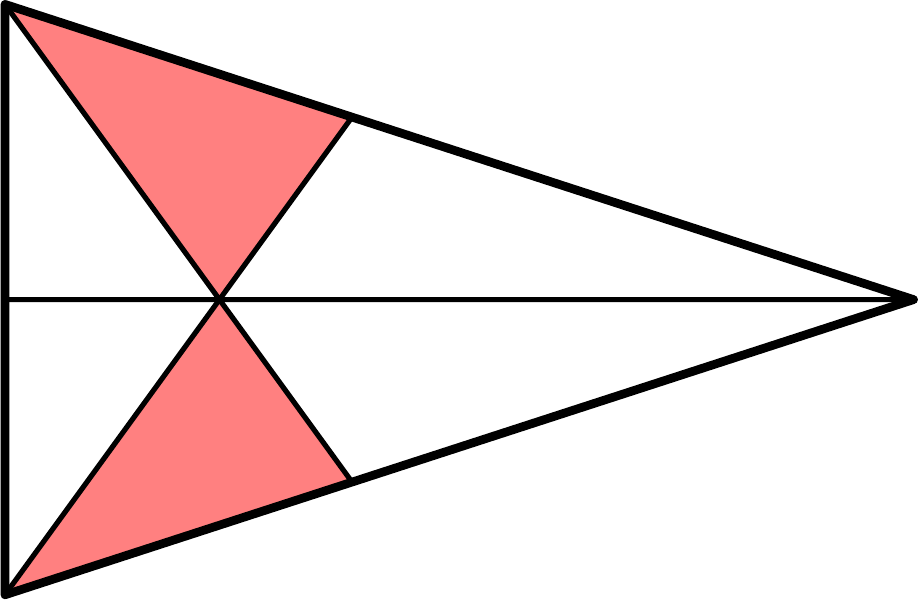}
\caption{The triangle $({\pi\over 5},{\pi\over 5},{2\pi\over5})^*$ subdivided using the incenter (the shaded triangles are similar to the original).}
\label{fig:eigen}
\end{figure}

Very little is known about the limiting distribution.  Experimentally, it appears that the densest point on the limiting distribution (i.e., the darkest region in Figure~\ref{fig:bb12}) corresponds to the triangle $({\pi\over 5},{\pi\over 5},{2\pi\over5})^*$.  This is likely because this triangle corresponds to an eigenvector of eigenvalue $1$ of two of the $M_i$.  In other words under subdivision using the incenter this triangle has two of its daughters which are similar to it (see Figure~\ref{fig:eigen}).  No other triangle has this property, and the triangle $({2\pi\over 9},{3\pi\over 9},{4\pi\over9})^*$ is the only other triangle with one of its daughters similar to itself, but this triangle does not appear to play a significant role in the distribution.

If we were to draw the histogram for $n=20$ or $n=50$ and compare it to Figure~\ref{fig:bb12} we would see almost no perceptible difference between them.  This is because, as we noted above, the convergence to the limiting distribution is exponential.  Or put another way, if we look at what happens when we map a triangle under $n$ applications of the $M_i$ in $\mathbf{P}$ then knowing what the last few steps that we applied gives us a good handle on where we are in $\mathbf{P}$.  This is essentially the heart of the proof of Theorem~\ref{thm:incenter}.

For example, the region $M_1M_2M_3\mathbf{P}$ corresponds to a triangle with vertices in $\mathbf{P}$ at $({\pi\over4},{\pi\over4},{\pi\over2})^*$, $({\pi\over8},{\pi\over4},{5\pi\over8})^*$ and $({\pi\over8},{\pi\over8},{3\pi\over4})^*$.  In particular, looking at all the daughters for $n$ large at least $1/6^3$ of the daughters must lie in this subregion of $\mathbf{P}$ (i.e., $1/6^3$ of the possible products of the $M_i$ will have $M_1M_2M_3$ as the leading term).  Since points inside this subregion of $\mathbf{P}$ must have minimum angle at least $\pi/8$ then we have that at least $1/6^3$ of the daughters in the $n$th generation must have minimum angle at least $\pi/8$.

\begin{figure}[ht]
\centering
\hfil\subfloat[Second generation]{\includegraphics[width=1.66in]{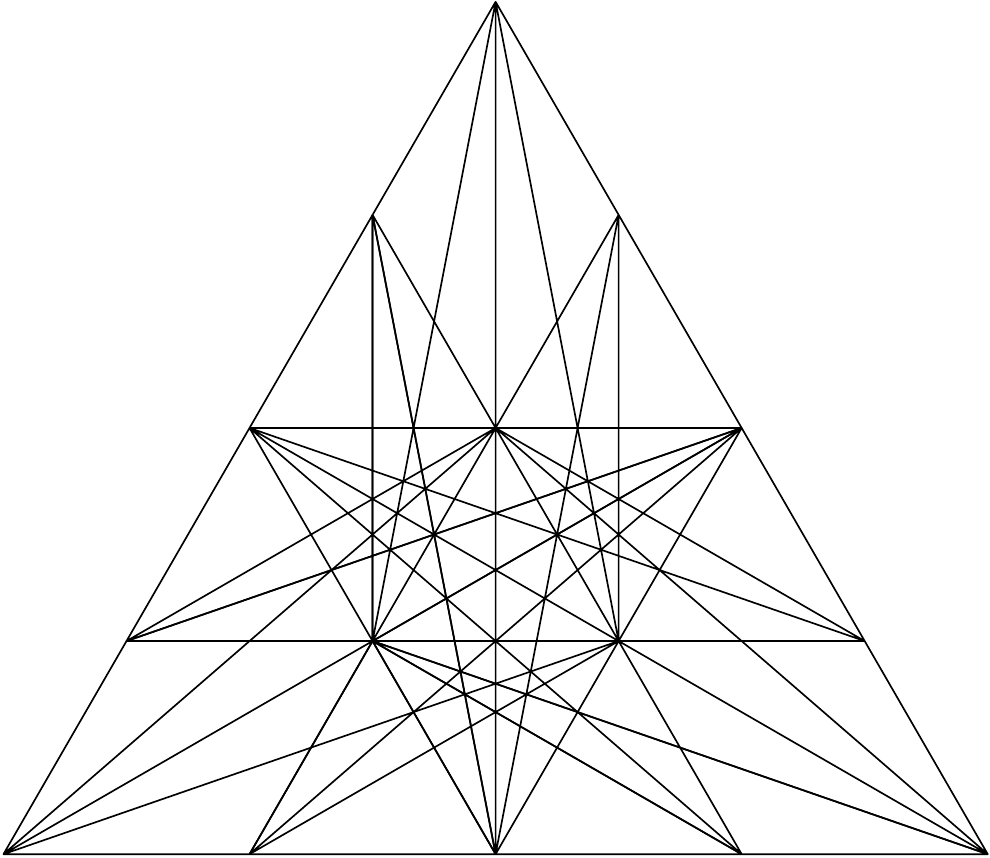}\label{fig:tta}}
\hfil\hfil\subfloat[Third generation]{\includegraphics[width=1.66in]{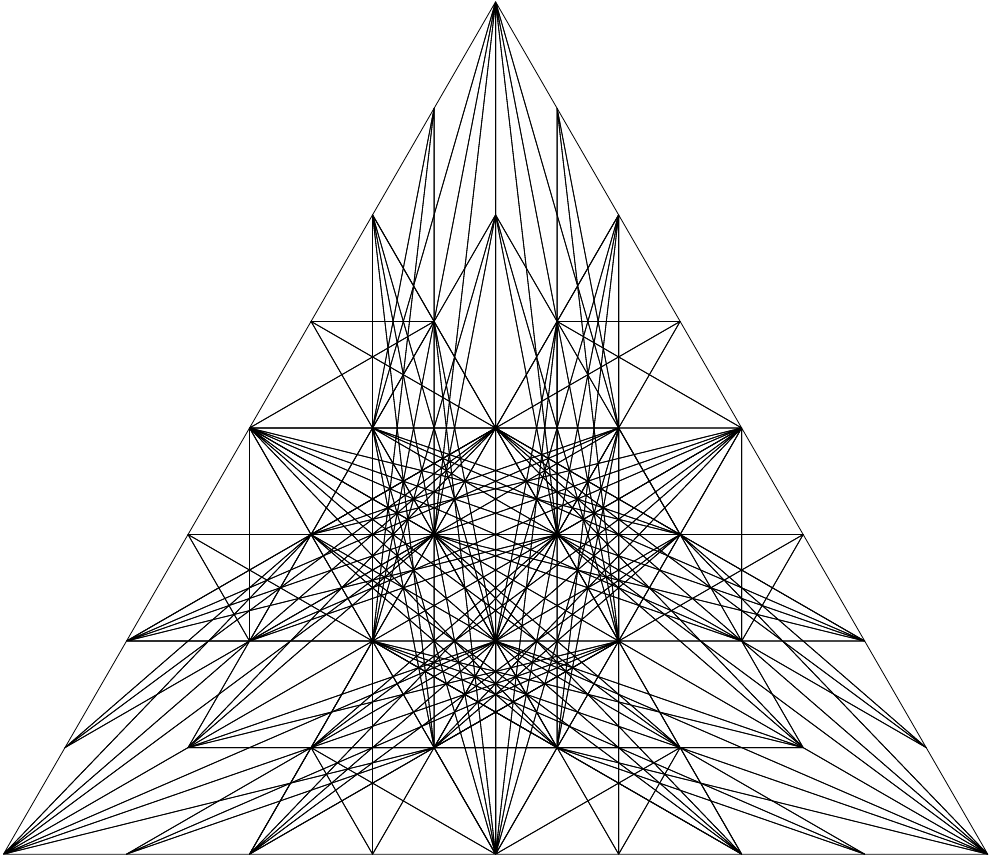}\label{fig:ttb}}
\hfil\hfil\subfloat[Fourth generation.]{\includegraphics[width=1.66in]{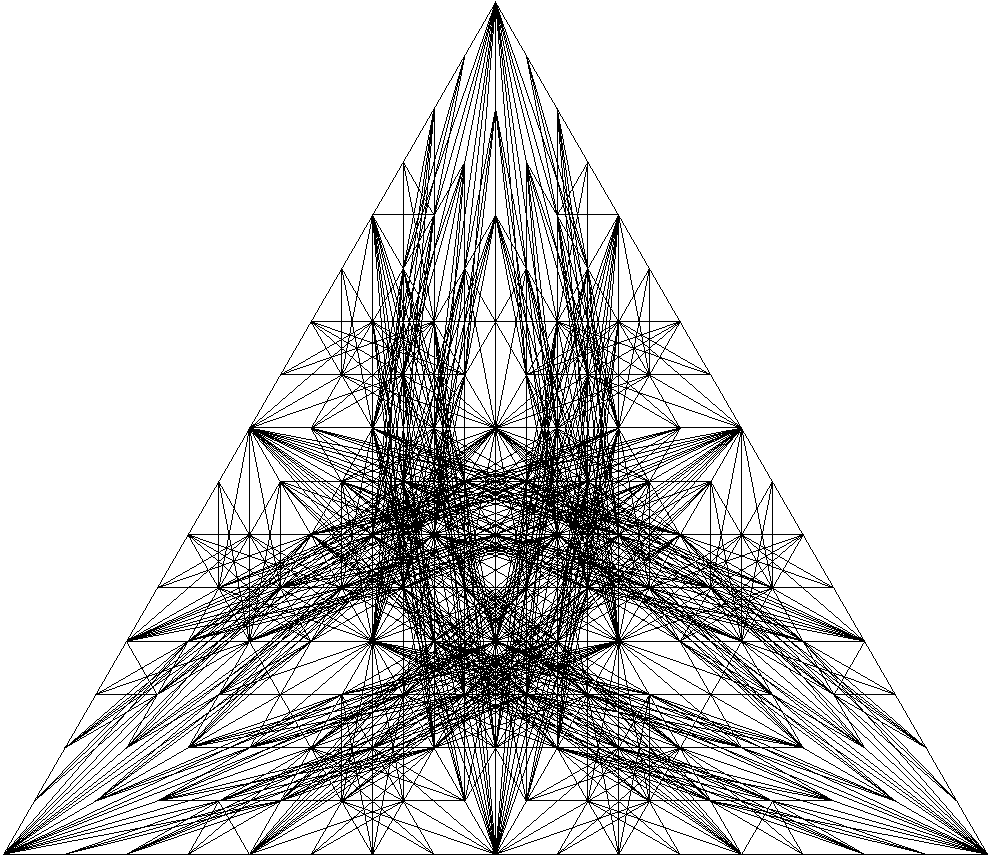}\label{fig:ttc}}
\hfil
\caption{Resulting images of $\mathbf{P}$ for different generations.}
\label{fig:trigen}
\end{figure}

Of course, by looking at larger products and looking over more prefixes we can say a lot more about what happens with the minimum angle.  For example, in Figure~\ref{fig:tta},~\ref{fig:ttb} and~\ref{fig:ttc} we have plotted all of the triangles that are the boundary of $P$ under all possible second ($M_iM_j\mathbf{P}$), third ($M_iM_jM_k\mathbf{P}$) and fourth ($M_iM_jM_kM_\ell\mathbf{P}$) maps of $\mathbf{P}$ into itself.  (The number of such triangles is large so it is difficult to pick out the individual triangular regions.  It also is interesting to note the similarity with the histogram in Figure~\ref{fig:bb12}.)

\begin{figure}[ht]
\centering
\includegraphics[scale=0.8]{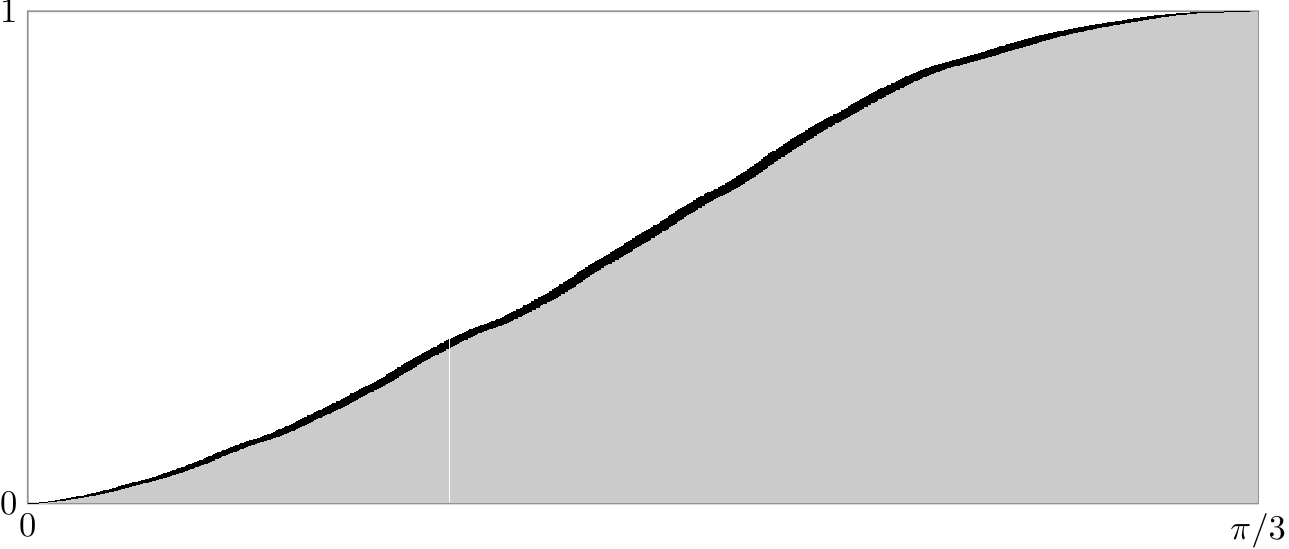}
\caption{Upper and lower bounds for the CDF for the smallest angle in the limiting distribution.}
\label{fig:CDF}
\end{figure}

We can now bound the limiting cumulative distribution function (CDF) for the smallest angle in the limiting distribution of triangles.  This is done by considering all resulting $6^n$ images of $\mathbf{P}$ in the $n$th generation.  Then for any angle $\theta$, a lower bound for the number of triangles with minimum degree $\theta$ (or less) is found by counting the number of the images of $\mathbf{P}$ which have largest minimum degree at most $\theta$.  Similarly, an upper bound for the number of triangles with minimum degree at most $\theta$ (or less) is found by counting the number of the images of $\mathbf{P}$ which contain a triangle with minimum degree at most $\theta$.  Doing this for $n=11$ gives Figure~\ref{fig:CDF}.  (By comparison the limiting CDF under subdivision using centroids is the constant function $1$, showing that these two methods of subdividing are fundamentally different.)

\section{Conclusion}\label{sec:conclusion}
We have seen that like the centroid, when doing subdivision using the incenter the resulting triangles are dense in the space of all triangles.  However, unlike the centroid the smallest angle in a typical triangle does not tend to $0$.  This is important since  certain methods can fail when the subdivision creates a large number of triangles with minimal angles going to $0$ as $n$ gets large (see \cite{bern,rivara,rosenberg}).

One interesting question is to understand the limiting distribution for repeated subdivision using the incenter, an approximation of which is shown in Figure~\ref{fig:bb12}.

\begin{figure}[ht]
\centering
\hfil\subfloat[Gergonne point]{\includegraphics[width=2.25in]{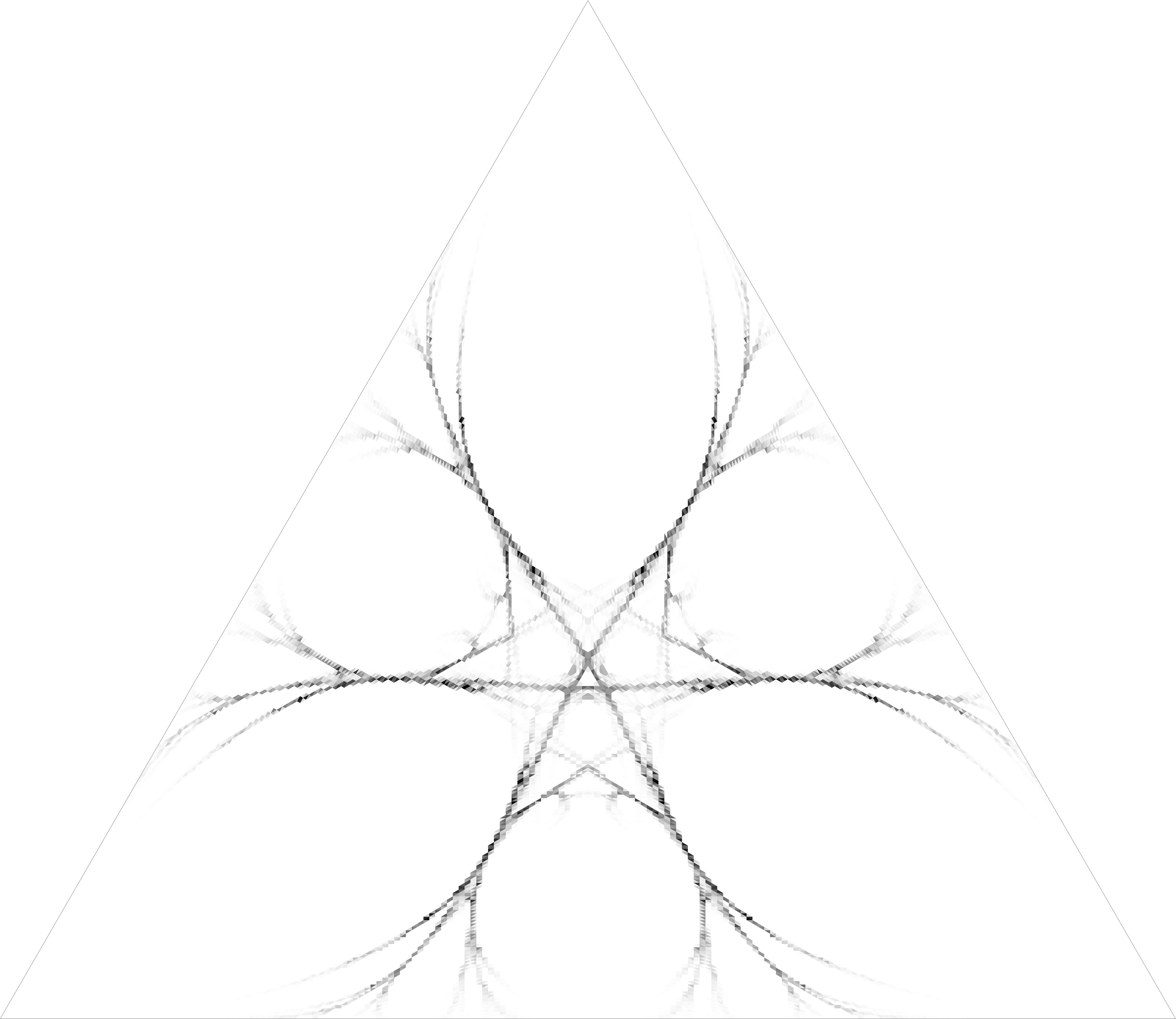}\label{fig:gerg}}
\hfil\subfloat[Lemoine point]{\includegraphics[width=2.25in]{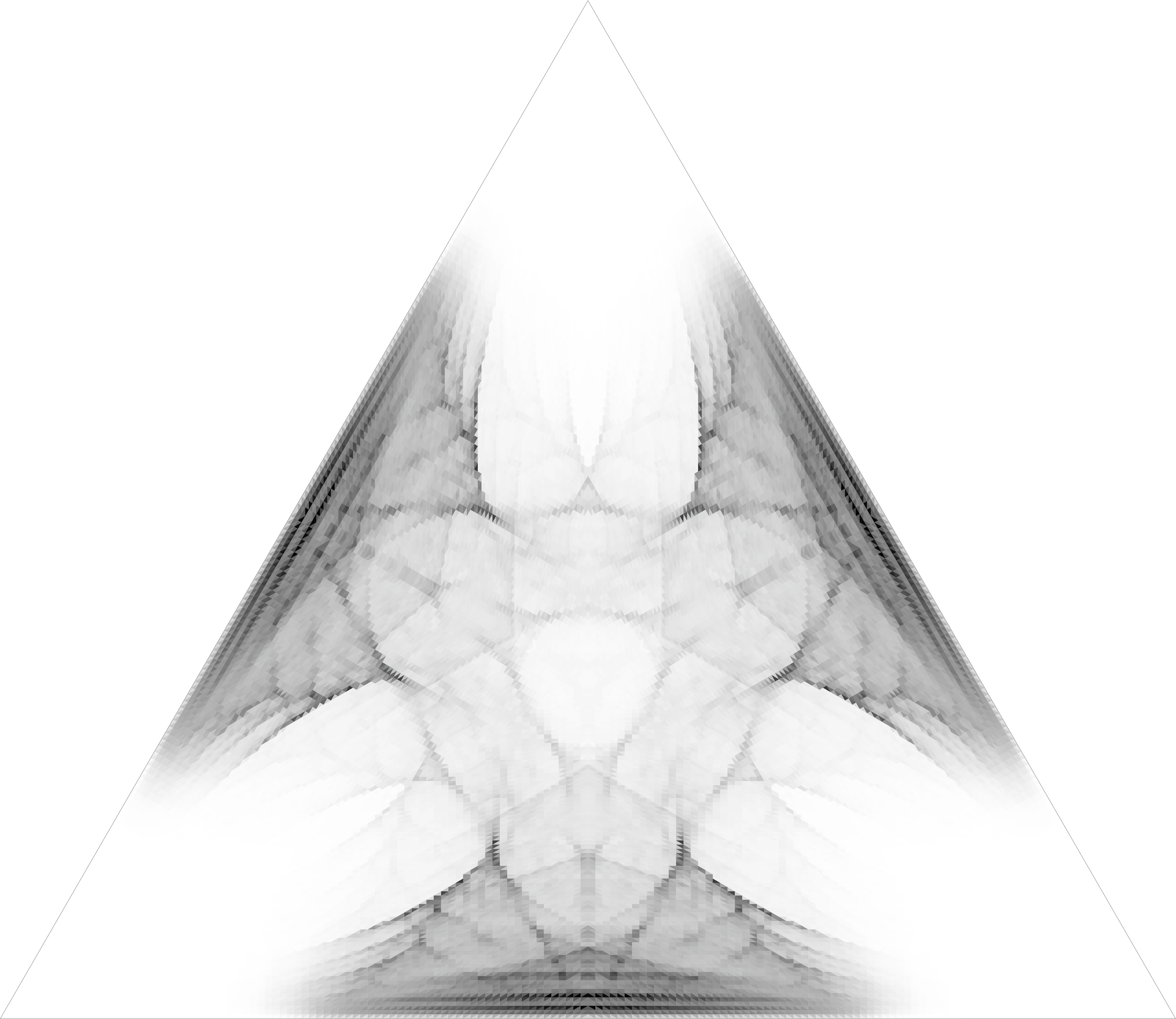}\label{fig:lem}}\hfil
\caption{Histograms for the Gergonne and Lemoine points.}
\end{figure}

One can also consider what  happens for subdivision using other interior points.  For example, the Gergonne point is found by taking the inscribed circle in the triangle and connecting a vertex to the point of tangency on the opposite edge; these three lines intersect at the Gergonne point.  When using the Gergonne point to subdivide it is known \cite{butler} that the triangles are not dense in the space of all triangles.   In Figure~\ref{fig:gerg} we have given a histogram of $\mathbf{P}$ for the tenth generation of subdividing using the Gergonne point (notice the large white spaces where there are no triangles).  

An interesting point for which little is known about what happens after repeated subdivision is the Lemoine point, which is found by taking the lines from a vertex through the median of the opposite edge and then flipping them across the angle bisectors; these three lines intersect at the Lemoine point.  In Figure~\ref{fig:lem} we have given a histogram of $\mathbf{P}$ for the eleventh generation of subdividing using the Lemoine point.  It is currently unknown if this method of subdivision is dense in the space of all triangles and what the limiting behavior is (there is some experimental evidence that the triangles become flat, but the convergence seems to be relatively slow).

More information about the Gergonne and Lemoine points, as well as a large number of other interesting points available to investigate, can be found online (see \cite{kimberling}).  More information about what happens under repeated subdivision using a central point can be found in \cite{iterated}.

\end{document}